\documentclass[a4paper,12pt,draft]{article}
\usepackage{mathtext}
\usepackage[T1,T2A]{fontenc}
\usepackage[cp1251]{inputenc}
\usepackage[english]{babel}
\usepackage{amsmath}
\usepackage{amsfonts}
\usepackage{amssymb}
\usepackage{mathrsfs}
\usepackage{amsthm}
\newcommand{\eps}{\varepsilon}
\newtheorem{cor}{Corollary}
\DeclareMathOperator{\essinf}{essinf}
\DeclareMathOperator{\esssup}{esssup}
\DeclareMathOperator{\spec}{spec}
\DeclareMathOperator{\supp}{supp}
\DeclareMathOperator{\Ann}{Ann}

\theoremstyle{plain}\newtheorem{T1}{Theorem}
\theoremstyle{plain}\newtheorem{T2}[T1]{Theorem}
\theoremstyle{plain}\newtheorem{L1}{Lemma}
\theoremstyle{plain}

\begin{document}

\title{New correction theorems in the light of a weighted Littlewood--Paley--Rubio de Francia inequality}
\author{D.M.Stolyarov}
\maketitle

\begin{abstract}
We prove the following correction theorem: every function $f$ on the circumference $\mathbb{T}$ that is bounded by the $\alpha_1$-weight $w$ (this means that $Mw^2 \leq C w^2$) can be modified on a set $e$ with $\int\limits_{e} w \leq \eps$ so that its quadratic function built up from arbitary sequence of nonintersecting intervals in $\mathbb{Z}$ will not exceed $C \log \frac{1}{\eps} w$.
\end{abstract}

\section{Introduction}

Correction theorems assert that an arbitary measurable function can be modifed on a set of a small measure up to a function with some good properties. Seemingly, the first and the most popular theorem of this type is the classical Lusin theorem about correction up to a continuous function. The next step was D.E.Men'shov's theorem \cite{Men} about correction of a bounded measurable function on a set whose measure does not exceed $\eps$ up to a function whose partial Fourie sums do not exceed $\frac{C}{\eps}$ uniformly and whose Fourier series converges uniformly. In 1979 S.V.Kislyakov in \cite{Ki4} sharpened the estimate of partial sums up to $C \log \frac{1}{\eps}$ and invented a general method of proving correction theorems. In the paper \cite{Ki2} the estimate was refined. 

We consider the circumference equiped with some weight $a(x), x \in \mathbb{T},$ as our general measure space. Similar statements hold for the line, but a slight change of technical details is needed. Now we turn to formal presentation, but first we need some definitions. 

\subsection{Definitions}

First, the Muckenhoupt conditions $A_p, 1 \leq p \leq \infty,$ will play a significant role in what follows. For every number $p$, $1 < p < \infty$, this condition can be written as follows:

 \begin{equation}
 A_p: \sup_{I} \big(\frac{1}{|I|}\int\limits_{I} w \big) \big(\frac{1}{|I|}\int\limits_{I}w^{-\frac{1}{p-1}} \big)^{p-1} < \infty.
 \end{equation}  
 
 If $p=1$, this condition will turn into $Mw \leq C w$ with some constant $C$, where $M$ is the Hardy-Littlewood maximal operator. We say that the weight satisfies the condition $A_{\infty}$ if it satisfies $A_p$ for some $p$. The theory of weights that obey such conditions can be found in the book \cite{St}. We also need a more sophisticated condition, $\alpha_p, 1 \leq p \leq 2$, which was introduced in the paper \cite{Ki1}. Specially, for $1 < p < 2$ a weight $w$ satisfies $\alpha_p$ if 
 
 \begin{equation}
 \alpha_p: \sup_{I} \big(\frac{1}{|I|}\int\limits_{I} w^{-\frac{1}{p-1}}\big)^{p-1} \big(\frac{1}{|I|}\int\limits_{I}w^{\frac{2}{2-p}} \big)^{\frac{2-p}{2}} < \infty.
 \end{equation}  
 
The supremum is taken over the set of all arcs of the circumference. It is easy to see that this condition is equalent to $w^{-\frac{1}{p-1}} \in A_{\frac{p'}{2}}$, or similary, $w^{\frac{2}{2 - p}} \in A_{\frac{p}{2-p}}$, where $p'$ is the exponent conjugate to $p$. The condition $\alpha_p$ can be extended to the border cases of $p=1$ or $p=2$ by passing to the limit, and $\alpha_1$ and $\alpha_2$ read as follows:

 \begin{equation}
 \alpha_1: w^2 \in A_1; \quad \alpha_2: w^{-1} \in A_1.
 \end{equation}
We see that $A_p$ follows from $\alpha_p$. Indeed, if $w \in \alpha_p$, then $w^{-\frac{1}{p-1}} \in A_{\frac{p'}{2}}$. Since the Muckenhoupt classes increase as the index increases, we also have $w^{-\frac{1}{p-1}} \in A_{p'}$, which is equalent to the inclusion $w \in A_p$. The case of $p=1$ can be obtained with the help of the Cauchy-Schwarz inequality. 

Second, we need the concept of the quadratic function $\sigma$. Let $\Delta_j, j \in \mathbb{N},$ be a family of disjoint segments in $\mathbb{Z}$. For each of them we introduce the corresponding Fourier multiplier $M_{\Delta_k}$ with the help of the following formula:

\begin{equation}
 M_{\Delta_k}(f) = (\chi_{\Delta_k}(\xi)\hat{f}(\xi)) \check{}.
 \end{equation}  
 
The formula is consistent even if $f$ is a distribution, and a fortiori if $f \in L^1(\mathbb{T})$. Now we can form the quadratic function:

 \begin{equation}
 \sigma f (x) = \big(\sum \limits_{k \in \mathbb{N}}|(M_{\Delta_k}f)(x)|^2 \big)^{\frac{1}{2}}.
 \end{equation}  
 
It was proved in the fundamental paper \cite{RdeF} that $||\sigma f||_{L^p(w)} \leq C ||f||_{L^p(w)}$ for all $w \in A_{\frac{p}{2}}$ when $2<p<\infty$. We will mostly use the result of the paper \cite{Ki1}, which is somehow dual to the previous one and can be written as follows. 

{\itshape
Suppose $1<p<2$, $0<r<p$ and $w \in \alpha_p$};{\itshape let $f_k$ be a sequence of summable functions such that $\supp \hat{f}_k \subset \Delta_k$. Then the following inequality holds}:{\itshape

$$||\sum\limits_k f_k||_{L^r(w)} \leq B_r ||\big( \sum\limits_k |f_k|^2 \big)^{\frac{1}{2}}||_{L^r(w)},$$
where $B_r$ does not depend on $\{f_k\}$ }({\itshape it only depends on the constant in the $\alpha_p$-estimate of $w$ and $r$}).

This theorem still holds for the cases of $p=1, p=2$, they are Corollaries $1,2$ in \cite{Ki1}.

\subsection{Statement of the main result}

Now we are ready to formulate the main result.

\begin{T1}
Suppose a weight $a$ satisfies the $A_{\infty}$ condition, and a weight $w$ satisfies the $\alpha_1$ condition. Let $f$ be a measurable function such that $|f| \leq w$. Then for every $\eps$, $0<\eps<1,$ there exists a function $g$ such that
$|g|+|f-g| = |f|$ and the following inequalities hold\rm{:}

$1)$ $\int\limits_{\{f \ne g\}} a \leq \eps \int\limits_{\mathbb{T}}|\frac{f}{w}|a,$
  
$2)$ $\sigma g \leq C(a,w) (1 + |\log(\eps)|) w.$ 

\end{T1}

From the condition $|g|+|f-g| = |f|$, it follows that correction is done by multiplying the initial function by some real nonnegative function $\phi$ whose values do not exceed one. The first inequality estimates the measure of the set where we correct the function. For example, if we take $w=a$, the measure of this set will be estimated by the Lebesque $L^1$-norm of $f$. The second inequality gives a pointwise estimate of the quadratic function in terms of the weight $w$. For example, we can try to make the weight $w$ sufficiently small (but it should still be separated from zero, otherwise it will not satisfy the $\alpha_1$-condition) on some set. In the fourth section, we will discuss special consequences of Theorem $1$ in detail. This theorem looks like Theorem $2'$ in \cite{AnKi}; in a way it is a generalization of that theorem, because only special sequences of disjoint intervals were involved there, but we have an arbitary one. On the other hand, we should pay for such a generality and the price is the condition on weight $w$, $A_1$ turned into $\alpha_1$ which, as we know, is stronger. We also mention that in \cite{AnKi} the logarithm in the estimate was squared, in our formula it is not. 

We are going to prove this theorem via the general method of obtaining correction theorems, which was described in \cite{Ki2}. We will need a weak $(1,1)$-type inequality for some operator, it will stated in the next subsection.

\subsection{An inequality}

Suppose $\mu$ is a measure, then we will denote by $L^p(l^2,\mu)$ the space of functions with values in $l^2$ that are summable in the $p$-th power with respect to $\mu$. Consider the operator $T$ defined on the set of finite sequences of trigonometric polynomials by the following formula:

\begin{equation}
 T(\{f_j\}) = \sum_j M_{\Delta_j}f_j.
\end{equation}

We are going to use another operator $T_u$, intertwined with $T$ with the help of multiplication by $u$. Specifically

 \begin{equation}
 T_u(\{f_j\}) = u^{-1}T(\{uf_j\}). 
 \end{equation}
 
Now we can formulate the second result of this paper.

\begin{T2}
Suppose a weight $a$ satisfies the classical condition $A_{\infty}$, a weight $w$ satisfies condition $\alpha_1$, see the first section. Let $u = \frac{a}{w}$. Then the operator $T_u$ defined by formulas $(6),(7)$ is continuous from $L^1(l^2,a)$ to $L^{1,\infty}(a)$.
\end{T2}

Strictly speaking, the statement needs further explanations, because we have defined the operator $T$ on the set of trigonometric polynomials, but now we apply it to some other functions. But as usual, it will be seen from the proof that everything is consistent. This theorem looks like theorem $4$ in \cite{AnKi}, but in that theorem instead of an operator $T$ there was a singular integral operator. $T$ is not an operator of that type, though it can be obtained as a composition of singular integral operators with somewhat nonstandard conditions on the kernel. Of course, the weak $(1,1)$-type conditions could have been destroyed under composition, fortunately, this does not happen.

Some words are in order about the operation of multiplication or division by $u$. An isometry between $L^1(a)$ and $L^1(w)$ is established in this way, both in the case of scalar-valued and $l^2$-valued functions. However, the operation fails to establish an isometry between the corresponding weighted Lorentz spaces $L^{1,\infty}$.

Before we turn to the proofs, we should make three small remarks. First, during the proof we assume that all sequences of functions are finite. It will allow us not to think about various technical convergence questions. The general case can be obtained by passing to the limit. Second, we assume all segments $\Delta_k$ to be contained in $\mathbb{Z}_{+}$, which will allow us to formulate the theorem in terms of analytic Hardy classes $H^p_{A}$. The general case can be obtained by adding the operators built up from the set of positive segments and the set of negative ones. If some of them contains zero, we can consider it by separately and then add it to the reminder.

As has already been mentioned, we are going to use Hardy classes. To be honest, the continuity of the operators mentioned above is related sooner to the properties of the Hardy spaces than of those of Lebesque spaces. We explain our notation. By $H^p_{A}(l^2,a)$ we denote the analytic Hardy class, which consists of all functions from the Smirnov class with values in $l^2$ whose boundary values are in $L^{p}(l^2,a(x) dx)$. We will often identify functions belonging to such classes with their boundary values. We will also use the $H_A^{1,\infty}(l^2,a)$ class, we think of it as of the closure of the set of finite sequences of analytic thigonometric polynomials in $L^{1,\infty}(l^2,a)$. This definition is nonstandard, usually $H_A^{1,\infty}(l^2,a)$ is defined as the intersection of $L^{1,\infty}(l^2,a)$ with the Smirnov class. For the $L^p$-norm, the definitions are equivalent when $p \geq 1$, but for the $L^{1,\infty}$-quasinorm the equivalence fails. For our purposes it will be more convenient to use the definition with trigonometric polynomials.

We are also going to use some interpolation technique to prove the second theorem. For the reader who is not familiar with it we can advice the book \cite{BL}; we also recall the notion of  $K$-closedness and its relation with interpolation. Let $(X_0, X_1)$ be a compatible couple of quasi-Banach spaces, and let $Y_0$ and $Y_1$ be closed subspaces of $X_0$ and $X_1$, respectively. The couple $(Y_0, Y_1)$ is said to be $K$-closed in $(X_0, X_1)$ if for every $y \in Y_0 + Y_1$ and a decomposition $y = x_0 + x_1, x_i \in X_i (i = 0,1)$ there exists another decomposition $y = y_1 + y_2$, where $y_i \in Y_i$ and $||y_i||_{Y_i} \leq C||x_i||_{X_i}$. It is easy to see that if $(Y_0, Y_1)$ is $K$-closed in $(X_0, X_1)$, then 

$$(Y_0, Y_1)_{\theta, q} = (Y_1 + Y_2) \cap (X_0, X_1)_{\theta, q}.$$
See \cite{Ki3} on the concept of $K$-closedness, its role in interpolation of Hardy spaces.

First, we will prove Theorem $2$, second, derive Theorem $1$ from it, and then discuss the meaning of these results. We turn to the proof.

\section{Proof of the second theorem}

The main ideas of the proof are similar to those employed in \cite{KiPar}. We are going to represent $T_u$ as a composition of a finite number of operators $\Tilde{T_u}$ with certain singular integral operators $\Tilde{T}$. We preface this by a lemma which is not related to the similarity transformation $\Tilde{T} \mapsto u^{-1}\Tilde{T}u$ itself, but plays a significant role in the proof.

\subsection{The weight mixing lemma}

\begin{L1}
Suppose $w \in \alpha_q,a \in A_{\infty}, 1< q <2$. Then there is $\delta$, $1 > \delta > 0$, such that for all $t$ in the interval $[1 - \delta, 1)$ there exists $r$ in $(1,2)$ such that the weight $w^t a^{1 - t}$ satisfies the condition $\alpha_r$. 
\end{L1}

\begin{proof}
We must estimate the quantity:

$$\big(\frac{1}{|I|} \int\limits_{I} w^{-\frac{t}{r-1}}a^{\frac{t-1}{r-1}} \big)^{r-1} \big(\frac{1}{|I|} \int\limits_{I} w^{\frac{2t}{2-r}} a^{\frac{2-2t}{2-r}} \big)^{\frac{2-r}{2}},$$
where $r$ is to be chosen.

Since $a \in A_{\infty}$, there exists $p$ such that $a \in A_p$. Therefore, by the Jones factorization theorem (see \cite{St}, Chapter $5$), there exist $a_1,a_2 \in A_1$ such that $a = a_1 {a_2}^{1-p}$. Now we substitute this new representation for $a$ in the formula and rewrite it in a bit different manner:

\begin{align*} \big(\frac{1}{|I|} \int\limits_{I} {\frac{1}{w^{\frac{t}{r-1}}}}{\frac{1}{a_1^{\frac{1-t}{r-1}}}} {a_2}^{\frac{(1-t)(p-1)}{r-1}} \big)^{r-1} \big(\frac{1}{|I|} \int\limits_{I} w^{\frac{2t}{2-r}} {a_1}^{\frac{2-2t}{2-r}} {\frac{1}{a_2^{\frac{(2-2t)(p-1)}{2-r}}}} \big)^{\frac{2-r}{2}} \leq \\
 \frac{1}{\essinf_{I} {a_1}^{1-t}} \frac{1}{\essinf_{I} {a_2}^{(p-1)(1-t)}} \big(\frac{1}{|I|} \int\limits_{I} \frac{{a_2}^{\frac{(p-1)(1-t)}{r-1}}}{w^{\frac{t}{r-1}}} \big)^{r-1} \big(\frac{1}{|I|} \int\limits_{I} w^{\frac{2t}{2-r}} {a_1}^{\frac{2 -2t}{2-r}} \big)^{\frac{2-r}{2}}.
\end{align*}

We estimate each integral separately. We use the standard H{\"o}lder inequality with the exponents $\frac{r-1}{(1-t)(p-1)}, \frac{r-1}{r-1 - (1-t)(p-1)}$ for the first integral:

\begin{align*}
\big(\frac{1}{|I|} \int\limits_{I} \frac{{a_2}^{\frac{(p-1)(1-t)}{r-1}}}{w^{\frac{t}{r-1}}} \big)^{r-1} \leq \big(\frac{1}{|I|} \int\limits_{I} a_2 \big)^{(1-t)(p-1)} \big(\frac{1}{|I|} \int\limits_{I} \frac{1}{w^{\frac{t}{r-1 - (1-t)(p-1)}}} \big)^{r-1 - (1-t)(p-1)}.
\end{align*}

To use the H{\"o}lder inequality, we need that $\frac{r-1}{(1-t)(p-1)} \geq 1$. We remember this condition. Now we use the H{\"o}lder inequality with the exponents $\frac{2 -r}{2 - 2t}, \frac{2-r}{2t-r}$ for the second integral:

\begin{align*}
\big(\frac{1}{|I|} \int\limits_{I} w^{\frac{2t}{2-r}} {a_1}^{\frac{2 -2t}{2-r}} \big)^{\frac{2-r}{2}} \leq \big(\frac{1}{|I|} \int\limits_{I} w^{\frac{2t}{2t-r}} \big)^{\frac{2t-r}{2}} \big(\frac{1}{|I|} \int\limits_{I} a_1 \big)^{1-t}.
\end{align*}

Also, here we should require that $\frac{2 -r}{2 - 2t} \geq 1$. 
 
So, we see that the contribution of the weights $a_1$ and $a_2$ to the formula can be estimated by their $A_1$-constants in the powers $1-t$ and $(p-1)(1-t)$, respectively. It only remains to estimate the following:

\begin{align*}
\big(\frac{1}{|I|} \int\limits_{I} \frac{1}{w^{\frac{t}{r-1 - (1-t)(p-1)}}} \big)^{r-1 - (1-t)(p-1)} \big(\frac{1}{|I|} \int\limits_{I} w^{\frac{2t}{2t-r}} \big)^{\frac{2t-r}{2}}.
\end{align*}  

Set $r = tq$. If $t$ is sufficiently close to $1$, then $r$ is also in $(1,2)$, therefore, this specification for $r$ is permitted. So we can rewrite the above expression in the form

\begin{align*}
\big(\frac{1}{|I|} \int\limits_{I} \frac{1}{w^{\frac{t}{tq-1 - (1-t)(p-1)}}} \big)^{tq-1 - (1-t)(p-1)} \big(\frac{1}{|I|} \int\limits_{I} w^{\frac{2}{2-q}} \big)^{\frac{2-q}{2} t}.
\end{align*} 

Now, since $w \in \alpha_q$, we can conclude that $w^{ - \frac{1}{q-1}} \in A_{\frac{q'}{2}}$. Therefore the reverse H{\"o}lder inequality is valid for the weight $w^{-\frac{1}{q-1}}$ for some $s$. Also we note that $\lim_{t \rightarrow 1}\frac{t}{tq-1 - (1-t)(p-1)} = \frac{1}{q-1}$. What is more, this value is greater than $\frac{1}{q-1}$. So, for all $t<1$ in some neighbourhood of $1$ we can write the following estimate:

\begin{align*}
\big(\frac{1}{|I|} \int\limits_{I} \frac{1}{w^{\frac{t}{tq-1 - (1-t)(p-1)}}} \big)^{tq-1 - (1-t)(p-1)} \leq c \big(\frac{1}{|I|} \int\limits_{I} w^{\frac{-1}{q-1}} \big)^{t(q-1)}.
\end{align*} 

As a result, after substituting this estimate in the previous one, we get exactly the $\alpha_q$ condition for $w$, raised to the power $t$. We also have to check two remembered inequalities. First, the number $\frac{r-1}{(1-t)(p-1)} = \frac{tq-1}{(1-t)(p-1)}$ must be less than one. As $t \rightarrow 1$, this value goes to infinity and, eventually, will exceed one. Second, the number $\frac{2 -r}{2 - 2t} = \frac{2 - tq}{2 -2t}$ must be greater than one. This value goes to infinity too, so finally the required inequality follows. We have proved the lemma.
\end{proof} 

Now we formulate Lemma $2$. Its statement is quite similar to the first lemma, the proofs' difference is only in that the reverse H{\"o}lder inequality is applied to another term is brackets. So, we omit it.

\begin{L1}
Suppose $w \in \alpha_q,a \in A_{\infty}, 1< q <2$. Then there exists $\delta$, $1 > \delta > 0$, such that for all $t$ in the interval $(1 ,1 + \delta]$ there exists $r \in (t,2)$ such that the weight $w^t a^{1 - t}$ satisfies the $\alpha_r$ condition. 
\end{L1}

It should be mentioned that the proof of Lemma $1$ shows that for $t$ sufficiently close to $1$, the weight $w^t a^{1-t}$ satisfies the condition $\alpha_{tq}$, the same is true for Lemma $2$. 

\begin{cor}
Both lemmas remain true for $q=1$.
\end{cor}

We can argue in the following manner: if the weight $w$ satisfies the $\alpha_1$ condition, it also satisfies the $\alpha_{1 + \delta}$ condition with some $\delta > 0$. This can be explained as follows: $w \in \alpha_1$ hence, $w^2 \in A_1$, then, by the reverse H{\"o}lder inequality, $w^{2+\eps} \in A_1$, therefore $w^{2+\eps} \in A_{1+\eps}$. But this exactly means that $w \in \alpha_{2 - \frac{2}{2+\eps}}$. So it satisfies the assumptions of lemmas.

\subsection{Auxiliary operators}
We now define two operators, $S$ and $R$, which came from \cite{KiPar}, they also played a significant role in \cite{Ki1}. We begin with the operator $S$. To define this operator, we assume all intervals $\Delta_j$ to be of length $2^l$, though each $l$ can occur several times. We name the set of those  $j$ whose length is equal to $2^k$ by $B_k$. Let $\xi \in (0,1)$ be some number, we think of it as of a number close to one, and let $\phi_k$ be trigonometric polynomials on the circumference that satisfy the following conditions borrowed from \cite{KiPar} (see conditions $(6),(7)$ respectively in that paper):

\begin{align}
\hat{\phi_m}(n) = 0 \quad {\rm for} \quad n \notin [0,2^m]; \quad |\hat{\phi_m}(n)| \leq 1; \\
|(\phi_m)^{(r)}(e^{i\sigma})| \leq C_{r,u} 2^{(r+1-u)m}\sigma^{-u} \quad {\rm for} \quad \sigma \in [-\pi,\pi], u>1, r \in \mathbb{Z}_+.
\end{align}

Here the differention is in the variable $\sigma$. These polynomials can also be chosen to satisfy an additional condition, namely: $\hat{\phi_m} = 1$ on $[(1 - \xi)2^{m-1}, (1 + \xi)2^{m-1}]$. The construction of such polynomials was discussed in \cite{KiPar} in detail. Let $\{h_j\} \in H^p_A(l^2,w)$, then we can define 

\begin{equation}
S(h)(x) = \sum_k \sum_{j \in B_k} e^{ia_j x}  (h_j * \phi_k)(x).
\end{equation}

The convolution is well defined, because $\phi_k$ lies in the Schwarz class and its Fourier transform has compact support, furthermore, we remind the reader that we have agreed to think that the set of the intervals is finite, therefore, the sum is finite too. We can also assume that $p \ne 1$, because, as it has been mentioned, if $w \in \alpha_1$, then $w \in \alpha_{1+ \eps}$. Therefore for all $p$, by Lemma $2$ in \cite{Ki1}, $S$ is continuous from $H_A^{r}(l^2,w)$ to $L^{r}(w)$ when $0 < r < p$ (in \cite{Ki1} everything happened on the line, but for the circumference the arguments are much the same). 

The operator $R$ is defined with the help of a family of trigonometric polynomials, namely, let $A$ be some number greater than one. Then there exist polynomials $\beta_j, j > 0$ such that the following two conditions are satisfied:

\begin{equation}
\hat{\beta_j} \geq 0, \sum_j \hat{\beta_j} = \chi_{\mathbb{Z} \backslash \{0\}}, \spec \beta_j \subset [A^{j-1},A^{j+1}],
\end{equation}

\begin{equation}
\forall r \in \mathbb{R}_+ \quad R \in \mathfrak{L}(H^r(l^2,w) \rightarrow H^r(l^2,w)),  R(\{f_k\}_k) = \{f_k * \beta_j\}_{k,j}.
\end{equation}
   
In fact, in this wording (but without weight, i.e., $w = 1$), this statement appeared in \cite{KiPar}, named Lemma $1$, in \cite{Ki1} it was redesigned slightly and adjusted to the line, it was named Lemma $3$ there. In \cite{Ki1} $r$ was not arbitary, but only from $0$ to $p$ if $w \in \alpha_p$. 

However, in the present paper we will deal not with the operators $S$ and $R$ directly, but with $S_u$ and $R_u$, consequently, we want to know that the latter two are continuous. We will not have to invent something new, the usual change of deinsity works. Nevertheless, we will have to interpolate over an unusual scale of spaces to achieve the continuity on the Lorentz class. 

\subsection{Continuity after a density change}

We notice that $f \leftrightarrow uf$ is an isometric bijection between $L^t(l^2,a)$, $L^t(a)$ and $L^t(l^2,w^ta^{1-t})$, $L^t(w^ta^{1-t})$ respectively ($t$ can be smaller than one). We introduce two auxiliary spaces, $E^t_u(l^2) = u^{-1}H_A^t(l^2,w^ta^{1-t})$ and $E^t_u = u^{-1}H_A^t(w^ta^{1-t})$. It is easy to see that $E^t_u, E^t_u(l^2)$ are subspaces of $L^t(a),L^t(l^2,a)$ respectively. The space $E^{1, \infty}_u$ can be defined as the closure in $L_{1,\infty}(a)$ of the set of analytic polynomials divided by $u$. Obviously, the operators $S_u, R_u$ specified by formula $(7)$ are well defined on a dense subsets of $E^t_u(l^2)$ (trigonometric polynomials divided by $u$). 

\begin{L1}
The operators $S_u$ and $R_u$ are continuous from $E^{1, \infty}_u(l^2)$ to $E^{1, \infty}_u$ and $E^{1, \infty}_u(l^2)$, respectively. 
\end{L1}

\begin{proof}
We will prove that the operators $S_u$ and $R_u$ are continuous from $E^{t}_u(l^2)$ to $E^{t}_u$ and $E^{t}_u(l^2)$, respectively, for $t$ smaller than one, but lying inside some neighbourhoodof it. We notice that by the definition of the spaces $E_u^t(l^2)$, it suffices to prove the continuity of $S$ and $R$ on the spaces $H_A^t(l^2,w^ta^{1-t})$. But by Lemma $1$, $w^ta^{1-t} \in \alpha_r$ for some $r \in (1,2)$. Therefore, we can use Lemmas $2,3$ from \cite{Ki1}. 

When $t$ is greater than one, but lies in some small neighbourhood of it, the operators $S_u$ and $R_u$ are also continuous from $E^{t}_u(l^2)$ to $E^{t}_u$ and $E^{t}_u(l^2)$ respectively, for similar reasons, one must merely use Lemma $2$ instead of Lemma $1$.

Finally, we extend the result to the case of $\textquoteleft t = (1, \infty)$' by interpolation. To do this, we notice that the couple $(E^t_u(l^2), E^s_u(l^2))$ is $K$-closed in $(L^t(l^2,a), L^s(l^2,a))$ if $t,s$ are sufficiently close to $1$. Actually, we have an isometry that sends $(E^t_u(l^2), E^s_u(l^2))$ onto $(H_A^t(l^2,w^t a^{1-t}), H_A^s(l^2,w^sa^{1-s}))$. Therefore, we must prove the $K$-closedness of the last-mentioned couple in $(L^t(l^2,w^t a^{1-t}), L^s(l^2,w^sa^{1-s}))$. To verify this, we use Theorem $3.3$ in \cite{Ki3}. It suffices to check that $L^t(l^2,w^t a^{1-t}), L^s(l^2,w^sa^{1-s})$ are $BMO$-regular lattices, which is true if $\log(w^t a^{1-t}) \in BMO, \log(w^sa^{1-s}) \in BMO$ by Corollary $3.1$ in the same paper. But these weights satisfy $\alpha_r$ for some $r$ by Lemmas $1,2$, therefore they satisfy the $A_r$ condition, and the logarithm of such a weight is always in $BMO$. For scalar-valued spaces everything is the same. 

From $K$-closedness it follows that $(E^t_u(l^2), E^s_u(l^2))_{\theta, \infty} = (E^t_u(l^2) + E^s_u(l^2)) \cap L^{1,\infty}(l^2,a) = u^{-1} N^+ \cap L^{1,\infty}(l^2,a)$ when $\frac{1}{1} = \frac{\theta}{t} + \frac{\theta}{s}$ ($N^+$ is the Smirnov class, the last identity can be obtained along the lines of a proof of $K$-closedness in \cite{Ki3}). But surely, the space on the right in this identity includes $E^{1, \infty}_u(l^2)$ as a closed subspace. For the  scalar-valued case everything is similar. Consequently, the operators $S_u, R_u$ are continuous on their domains as operators from $E^{1, \infty}_u(l^2)$ to $E^{1, \infty}_u(l^2)$ and $E^{1, \infty}_u$ respectively, so they can be extended by continuity to these spaces, which proves the lemma.   
\end{proof}

\subsection{The end of the proof}

We can reformulate Theorem $2$ as an inequality:

\begin{equation}
\big|\big|\sum\limits_{k} u^{-1}M_{\Delta_k}(uf_k)\big|\big|_{L^{1,\infty} (a)} \leq C \big|\big|\big(\sum\limits_k |f_k|^2 \big)^{\frac{1}{2}}\big|\big|_{L^1(a)}.
\end{equation}

We need the following inequality:

\begin{equation}
\big|\big|u^{-1}\big(\sum |M_{\Delta_k}(uf_k)|^2\big)^{\frac{1}{2}}\big|\big|_{L^{1,\infty}(a)} \leq C \big|\big|\big(\sum |f_k|^2\big)^{\frac{1}{2}}\big|\big|_{L^1(a)},
\end{equation}
which means that for the projection $P$ defined by the formula $P(\{f_k\}_k) = \{M_{\Delta_k}f_k\}_k$, the corresponding operator $P_u$ is continuous from $L^1(l^2,a)$ to $L^{1,\infty}(l^2,a)$. To prove $(14)$, we observe that $M_{\Delta_k}f_k = e^{ 2 \pi i a_k \cdot} P_+(e^{-2 \pi i a_k \cdot}f_k(\cdot)) - e^{2 \pi i b_k \cdot} P_+(e^{- 2 \pi i b_k \cdot}f_k(\cdot))$, where $\Delta_k = [a_k, b_k-1]$ and $P_+$ is the Reisz projection. Therefore we have represented $P$ in the form $P = U_1^{-1}P_+ U_1 - U_2^{-1}P_+ U_2$, where $U_1,U_2$ are operators of multiplication by something unimodular, consequently, isometric operators both on $L^1(l^2,a)$ and on $L^{1,\infty}(l^2,a)$. Obviously, a similar formula holds for $P_u$. We recall that $w$ is in the $\alpha_1$ class and, consequently, in $A_1$. Therefore, the operator $(P_+)_u$ is continuous from $L^1(l^2,a)$ to $L^{1,\infty}(l^2,a)$ by Theorem $4$ in \cite{AnKi}. We see that the desired continuity property holds for $P_u$, and the inequality is proved.  

So we must only prove the following: 

\begin{equation}
\big|\big|\sum u^{-1}M_{\Delta_k}(uf_k) \big|\big|_{L^{1,\infty} (a)} \leq C \big|\big| u^{-1}\big(\sum |M_{\Delta_k}(uf_k)|^2 \big)^{\frac{1}{2}} \big| \big|_{L^{1,\infty}(a)}.
\end{equation}

If we denote $g_k = u^{-1}M_{\Delta_k}(uf_k)$, we need to prove that the operator intertwined with the help of multiplication by $u$ with the operator $\{g_k\}_k \mapsto \sum\limits_k g_k$ acts from $E^{(1,\infty)}_u(l^2)$ to $E^{(1,\infty)}_u$ (recall that we deal with the case when all intervals $\Delta_k$ are contained in $\mathbb{Z}_+$) and is continuous. During the remaining part of the proof we will be busy with representing this operator as a composition of operators of type $R$ or $S$ for various sequences of intervals. Then the operator we require ultimately will be represented as a composition of operators like $R_u$ and $S_u$ and will be continuous. To be honest, this procedure was described in detail both in \cite{KiPar} and \cite{Ki1}. We repeat it here for completeness. 

Our first purpose is to make our sequence of intervals $\textquoteleft$more regular'. We will use a special $\textquoteleft$cutting' procedure. We move all our functions $f_j$ so that the left end of $\Delta_j$ goes to $1$, i.e., we introduce the functions $g_j(x) = e^{- i(a_j - 1)x}f_j(x)$, where $a_j$ is the left end of the interval $\Delta_j$. Then we apply the operator $R$ with $A = 2^{\frac{1}{10}}$ to $g$ and then return everything back. As a result, we get a set of functions $f_{jk}$ with more regular spectra, but with the same sum.

We first deal with those $j$ for which the length of $\Delta_j$ is at most eleven. We note that for any of such intervals the interval $9 \Delta_j$ does not intersect more than $99$ other intervals of the partition. Therefore, this set can be split into $100$ subsets in such a way that inside each set the spectra of the functions $f_j$ will be separated in a good way (this means that the segments $3 \Delta_j$ and $3 \Delta_{j'}$ do not intersect when $j \ne j'$, and $j,j'$ are contained in one subset of the partition). After that we apply the operator $S$ to this subsequence (to be more accurate, $S$ will be applied to $e^{-ia_jx}f_j(x)$) and get the part of the entire sum that is generated by these functions. As a result, we have got rid of small intervals.

Each of the remaining functions $f_j$ has been split into several functions $f_{jk}$. Define the set $A_m, m = [0,...,9], A_m = \{(j,k)| k \equiv_{10} m, f_{jk} \ne 0 \}$. That is, we have divided the functions $f_{jk}$ into ten groups in accordance with to the remainder after division $k$ by $10$. We note that if we disregard the functions that have the biggest and the biggest but one $k$ for each $j$, then, in each of the remaining groups, the spectra of functions are separated in a good way, because to the left and to the right of them there are at least two intervals from the other groups (we recall that the spectrum of $f_{jk}$ can intersect the spectrum of $f_{j(k-1)}$, but not of $f_{j(k-2)}$). Therefore we can again apply an operator of type $S$ to the set of functions from each group and get the part of the sum generated by these functions (again, we apply the $S$-type operator to functions ''shifted-to-zero''). As to  the functions we have omitted, we proceed similarly, but ''in the opposite direction''. Specially, we move the right ends of their intervals to $-1$, after that we apply an operator similar to $R$ but generated by antianalytic functions, and then return everything back. Then we do the the same procedure with partition in $10$ groups, the only difference is that now we do not need to avoid the small intervals, because the intervals of the partition were big enough and consequently, divided at least into four intervals. Therefore, there was something to the left from the last and the last but one interval, so the last intervals of the new partition are separated in a good way. As a result, their contribution to the entire sum can be rewritten in terms of application of an operator $S$ for some sequence. Now we see that we have represented the operator $\{g_k\}_k \mapsto \sum\limits_k g_k$ in the desired way, and, thus have proved Theorem $2$.

\section{The derivation of the first theorem from the second}

Consider the set $X$ of bounded measurable functions on the circumference for which $w^{-1}\sigma(wf) \in L^{\infty}$. We define the norm on this set by the formula:

\begin{equation}
||f||_X = \esssup\{|f(\cdot)|, w^{-1} (\sum\limits_k |M_{\Delta_k}(fw)|^2)^{\frac{1}{2}}\}.
\end{equation}

It is easy to see that $X$ is nonempty. Indeed, an $A_1$-weight is separated from zero, therefore, $X$ contains, for example, all functions that can be obtained from trigonometric polynomials via division by $w$. We still view $(\mathbb{T}, a(x) dx)$ as our main measure space. We are going to use a general theorem from \cite{Ki2} (we alter it a bit, to take into account to the fact that our space has finite measure):

{\itshape
Let a Banach space $X$ of $\mu$-measurable functions satisfy the following two conditions.

$A1$. The canonical embedding of $X$ into $L^1(\mu)$ is continuous and the unit ball of $X$ is weakly compact in $L^1(\mu)$.

$A2$. For every $g \in L^{\infty}$ the functional $\Phi_g$ on $X$ defined by the formula $\Phi_g(h) = \int gh d \mu$ satisfies the following inequality{\rm:}

$$ ||g||_{L^{1,\infty}(\mu)} \leq c ||\Phi_g||_{X^*},$$
where the constant does not depend on $g$.

Then for every function $F$ such that $||F||_{L^{\infty}} \leq 1$ and every $\eps, 0 < \eps < 1$, there a function $G$ such that $|G|+|F-G| = |F|$, $\mu(F \ne G) \leq \eps, ||G||_X \leq C(1 + \log \eps^{-1})$
} 

So, we have to check two conditions.

\subsection{The first condition}

$X$ embeds into $L^{\infty}$, consequently, it embeds into $L^1(a)$. We have to check that the unit ball of $X$ is compact in the weak topology of the space $L^1(a)$. We notice that the weak $L^1$-convergence on the ball of $L^{\infty}$ coincides with the weak* convergence in $L^{\infty}$ regarded as the dual of $L^1$; we will check the compactness in this last topology. We also see that $X$ is a subspace of $(L^{\infty} \bigoplus L^{\infty}(l^2))_{\infty}$. Indeed, we can define the embedding map $\alpha: X \rightarrow (L^{\infty} \bigoplus L^{\infty}(l^2))_{\infty}$ with the help of the following formula: $\alpha(h) = (h, \{w^{-1}M_{\Delta_k}(wh)\}_k)$. Therefore our ball $B_X$ becomes the image of the ball $\alpha B_X$ after the canonical projection $(L^{\infty} \bigoplus L^{\infty}(l^2))_{\infty}$ to the first coordinate. Thus we are to prove the compactness of $\alpha B_X$. This set is a subset of the ball of $(L^{\infty} \bigoplus L^{\infty}(l^2))_{\infty}$, but this space is conjugate to $(L^{1}(a) \bigoplus L^{1}(l^2,a))_{1}$, its ball is compact by the Alaoglu theorem, so we need only prove the closedness of $\alpha B_X$ viewed as a subset of the ball $(L^{\infty} \bigoplus L^{\infty}(l^2))_{\infty}$, in other words, we have to prove that if $f_n \rightarrow f$ weakly* in $L^{\infty}$ and $||f_n||_{X} \leq 1$, then $||f||_{X} \leq 1$, or, again the same, $||\alpha f||_{L^{\infty}(l_2)} \leq 1$. By definition, $\psi \mapsto M_{\Delta_k}(w \psi)$ is a continuous finite rank operator from $L^{\infty}$ to $C(\mathbb{T})$. Thus, $M_{\Delta_k}(wf_n) \rightarrow M_{\Delta_k}(wf)$ in $C(\mathbb{T})$, consequently, for every $N$ we have the estimate $w^{-1} (\sum\limits_{k = 1}^{N} |M_{\Delta_k}(wf)|)^{\frac{1}{2}} \leq 1$. Passing to the limit in $N$, we get the desired result. So, we have checked the first condition.

\subsection{The second condition}

We have to prove the estimate $||\Phi_g||_{X^*} \geq c||g||_{L^{1,\infty}(a)}$. $X$ is a closed subspace (as an isometric image) of $(L^{\infty} \bigoplus L^{\infty}(l^2))_{\infty}$. Now we prove that our functional is continuous in the topology induced on $X$ as on a subspace by the weak* topology of $(L^{\infty} \bigoplus L^{\infty}(l^2))_{\infty}$, viewed as the dual of $(L^{1}(a) \bigoplus L^{1}(l^2, a))_{1}$. As we know from the previous subsection, $X$ is a closed subspace in this topology, because its ball is closed (for example, we can use the Banach lemma see \cite{Ios}, addition to the $\S$ $5.4$). Let $\{h_k\}$ be a sequence in $X$. Let it converge to some $h$ in the above sense. Then we use the convergence of the first coordinates in $(L^{\infty} \bigoplus L^{\infty}(l^2))_{\infty}$ and see that $h_k \rightarrow h$ weakly* in $L^{\infty}$. But $g \in L^{\infty}$, $g \in L^1(a)$, therefore $\int h_k g \rightarrow \int h g$. The continuity is proved. As a result, this functional on $X$ can be identified canonicaly with an element of $(L^{1} \bigoplus L^{1}(l^2))_{1} / \Ann X$ so, we can choose a representative at which the norm is almost attained, i.e., a functional $\Tilde{\Phi}$ that extends $\Phi$ and satisfies $\Tilde{\Phi}((h,\{h_k\}_k)) = \int fha + \sum_k \int f_k h_k a$, where $(f, \{f_k\}_k) \in (L^{1}(a) \bigoplus L^{1}(l^2,a))_{1}$, $||(|f|^2 + \sum_k |f_k|^2)^{\frac{1}{2}}||_{L^1(a)} \leq ||\Phi||+ \eps$. Therefore, $\Phi(h) = \int f h a + \sum_k \int w^{-1}M_{\Delta_k}(w h) f_k a = \int f h a + \sum_k \int h u^{-1} \overline{M_{\Delta_k}( u \overline{f_k})}a$. Substituting trigonometric polynomials for $h$, we get $g = f + u^{-1} \sum_k \overline{M_{\Delta_k} (u\overline{f_k})}$. 

So we have to prove the following inequality:

\begin{equation}
\big| \big|f + u^{-1} \sum_k \overline{M_{\Delta_k} (u \overline{f_k})}\big|\big|_{L^{1, \infty}(a)} \leq c \big|\big|\big(|f|^2 + \sum_k |f_k|^2\big)^{\frac{1}{2}} \big|\big|_{L^1(a)}.
\end{equation}

Obviously, we can estimate $f$ and the remaining sum separately. But an estimate for $f$ is trivial with the constant one, and the inequality for the sum is preciesly Theorem $2$ in the form $(13)$. So we have checked the second condition too.

So the conditions of the quoted theorem are fulfilled. But the first theorem is absolutely similar to it, one only have to substitute $F = \frac{f}{w}$ instead of $F$. 

\section{Corollaries and a conjecture}

\subsection{Corollaries to Theorem $1$}

If we take $w = a = 1$ in the first theorem, we get the following statement. 

\begin{T1}
Let $f$ be a measurable function such that $|f| \leq 1$. Then for every $\eps, 0<\eps<1,$ there exists a function $g$ such that
$|g|+|f-g| = |f|$ and the following inequalities hold{\rm:}

$1)$ $\mu{\{f \ne g\}}  \leq \eps ||f||_{L^1(\mu)},$
  
$2)$ $\sigma g \leq C (1 + |\log(\eps)|)$.

\end{T1}

Now let $f$ be an arbitary measurable function from $L^{\infty}, |f| \leq 1$. We take $w = (Mf)^{\gamma}, 0 < \gamma < \frac{1}{2}$. Then $w^2 = (Mf)^{2\gamma} \in A_1$. So, we arrive at the following statement.	

\begin{T1}
Let $f$ be a measurable function such that $|f| \leq 1$, and let $a \in A_{\infty}, \gamma \in (0,\frac{1}{2})$. Then for each $\eps, 0<\eps<1$, there exists a function $g$ such that
$|g|+|f-g| = |f|$ and the following inequalities hold:

  $1)$ $\int\limits_{\{f \ne g\}}a(x) dx  \leq \eps \int f(x)^{1-\gamma}a(x) dx,$
  
  $2)$ $\sigma g \leq C (1 + |\log(\eps)|)(Mf)^{\gamma}$.

\end{T1}

\subsection{Theorem 1 on the line}

On the line the first theorem should be formulated in the following way.

\begin{T1}
Let $a$ satisfy $A_{\infty}$, let $w$ satisfy $\alpha_1$, and let $u = \frac{a}{w}$. Let $f$ be a measurable function with compact support such that $|f| \leq w$. Then for every $\eps, 0<\eps<1,$ there exists a function $g$ such that
$|g|+|f-g| = |f|$ and the following inequalities hold{\rm:}

  $1)$ $\int\limits_{\{f \ne g\}} a \leq \eps \int\limits_{\mathbb{T}}|\frac{f}{w}|a,$
  
  $2)$ $\sigma g \leq C(a,w) (1 + |\log(\eps)|) w$.

\end{T1}

The proof is absolutely the same, there is a small difference in a technical detail, because on the line a singular integral operators map $L^{\infty}$ only to $BMO$ and consequently, its values on $L^{\infty}$ can be defined only modulo constants. However, if we take the function we are going to correct from $L^{\infty} \cap L^p$, a singular integral operator will send it to $L^p$; for this purpose we impose the compact support condition.

\subsection{About the conditions on the weight $w$}

In this subsection we are concerned with several questions about weights in Theorems $1$ and $2$. The discussion is prefaced by the following lemma.

\begin{L1}
Suppose a weight $w$ satisfies the conditions $\alpha_p$ and $A_1$. Then $w \in \alpha_1$
\end{L1}

\begin{proof}

We raise the inequality that express the condition $A_1$ to some power $b$ (to be specified later) and multiply it by the inequality expressing $\alpha_p$. This results in the following inequality:

\begin{align}
\big( \frac{1}{|I|} \int\limits_I w^{-\frac{1}{p-1}}(x) dx\big)^{p-1} \big( \frac{1}{|I|} \int\limits_I w^{\frac{2}{2 - p}}(x) dx\big)^{\frac{2-p}{2}} \big( \frac{1}{|I|} \int\limits_I w(x) dx\big)^{b} \leq [w]_{\alpha_p}[w]_{A_1}\essinf_I(w^b).
\end{align}

Now we apply two H{\"o}lder inequalities. The first will be:

\begin{align}
\big( \frac{1}{|I|} \int\limits_I w^{-\frac{1}{p-1}}(x) dx\big)^{(p-1)a} \big( \frac{1}{|I|} \int\limits_I w^{\frac{2}{2 - p}}(x) dx\big)^{\frac{2-p}{2}} \geq \\
\geq \big( \frac{1}{|I|} \int\limits_I w^{-\frac{a}{c} + \frac{1}{c}}(x) dx\big)^{c} = \big( \frac{1}{|I|} \int\limits_I w^2(x) dx\big)^{c}.
\end{align}

The constants are: $c = (\frac{p-1}{p-1 + \frac{2-p}{2}} + 2)^{-1}, a = \frac{p-1}{p-1 + \frac{2-p}{2}} c$. Then the exponents of H{\"o}lder inequality are $\frac{c}{a(p-1)}$ and $\frac{2c}{2-p}$, they are conjugate indeed. What is more, we have $-\frac{a}{c} + \frac{1}{c} = 2$, which leads to the second identity.

The second H{\"o}lder inequality will look like this (to be preciese, this is a H{\"o}lder inequality raised to power):

\begin{align}
\big( \frac{1}{|I|} \int\limits_I w^{-\frac{1}{p-1}}(x) dx\big)^{(p-1)(1-a)} \big( \frac{1}{|I|} \int\limits_I w(x) dx\big)^{1-a} \geq 1.
\end{align}

Notice that, since $0<a<1$, we can take $b = 1-a$. Then, after multiplication of the first and the second H{\"o}lder inequalities, by using $(18)$ we get the following estimate:

$$\big( \frac{1}{|I|} \int\limits_I w^2(x) dx\big)^{c} \leq [w]_{\alpha_p}[w]_{A_1}\essinf_I(w^b).$$

We only have to check that $b = 2c$, then it will be exactly the $\alpha_1$-condition, raised to the power $c$. But this is so indeed, and the lemma is proved.
\end{proof}

It could have been thought that in Theorem $2$ one could require $w \in \alpha_p$ rather than $w \in \alpha_1$. Indeed, all the arguments remain valid, except for a small portion about the Reisz projection, between formulas $(14)$ and $(15)$. There we need $w \in A_1$, because the Reisz projection is discontinuous as an operator from $L^1(w)$ to $L^{1, \infty}(w)$ if $w \notin A_1$. This fact is well known, for example, see Proposition 5.4.7 in \cite{St}, which corresponds to a nearby situation. 

By Lemma $4$, these conditions together lead to $w \in \alpha_1$, so we cannot strengthen Theorem $2$ in such a way. 

\subsection{On interpolation}

Finally, we state an interesting conjecture, which was partly tackled during the proof of the second theorem and could have shorten it. Namely, we consider  the spaces $X^p$, which are obtained as the closure of the set of finite sequences of trigonometric polynoms that satisfy the conditions $\supp \hat{f_k} \in \Delta_k$, in the topology of $L^p(w)$. Then we suppose the following lemma to be true.

\begin{L1}
{\rm(}Conjecture{\rm)} Let $p_1<1<p_2$. Then the couple $(X^{p_1},X^{p_2})$ is K-closed in $(L^{p_1},L^{p_2})$.
\end{L1}

The work is supported by the Chebyshev Laboratory, grant 11.G34.31.0026 of the Government of Russian Federation.

The author is grateful to S.V.Kislyakov both for the posing the problem and for significant help in solving it.

\end{document}